# FILTRATION SHRINKAGE BY LEVEL-CROSSINGS OF A DIFFUSION

By A. Deniz Sezer

*York University*

We develop the mathematics of a filtration shrinkage model that has recently been considered in the credit risk modeling literature. Given a finite collection of points $x_1 < \cdots < x_N$ in $\mathbb{R}$, the region indicator function $R(x)$ assumes the value $i$ if $x \in (x_{i-1}, x_i]$. We take $\mathbb{F}$ to be the filtration generated by $(R(X_t))_{t \geq 0}$, where $X$ is a diffusion with infinitesimal generator $\mathcal{A}$. We prove a martingale representation theorem for $\mathbb{F}$ in terms of stochastic integrals with respect to $N$ random measures whose compensators have a simple form given in terms of certain Lévy measures $F_i^{j\pm}$, which are related to the differential equation $\mathcal{A}u = \lambda u$.

**1. Introduction.** Let $x_1, \ldots, x_N$ be a finite collection of points in $\mathbb{R}$, in increasing order. These levels separate $\mathbb{R}$ into $N + 1$ regions, namely, $(-\infty, x_1], (x_1, x_2], \ldots, (x_N, \infty)$. Let us put $R_i = (x_i, x_{i+1}], i = 0, \ldots, n-1$, $R_0 = (-\infty, x_1]$ and $R_N = (x_N, \infty)$. We define the region indicator function $R(x)$ as

$$R(x) = i \quad \text{if } x \in R_i.$$

Next, we consider a nonsingular diffusion $X$, with state space $I$, an interval of $\mathbb{R}$, and with infinitesimal generator $\mathcal{A}$ of the form

$$\mathcal{A} = \frac{1}{2} a(x) \frac{d^2}{dx^2} + b(x) \frac{d}{dx},$$

where $a(x)$ is strictly positive and continuous and $b(x)$ is locally integrable on $I$, acting on a domain of functions as described in [7, 10]. We assume that $I$ either is open or has absorbing boundary.

Assuming that the points $x_1, \ldots, x_N$ belong to the interior of $I$, we let $\mathbb{F}^0 = (\mathcal{F}_t^0)_{t \geq 0}$ be the filtration generated by $(R(X_t))_{t \geq 0}$. (We reserve the notation $\mathbb{F}$ for a slightly larger filtration, which is right continuous and complete.)

$\mathbb{F}^0$ is a sub-filtration of the natural filtration of $X$ and among the "filtration shrinkage" models that have recently been considered in the credit risk literature [5, 8, 15]. These models aim to represent the incomplete and dynamic information about a process of interest as a sub-filtration of the natural filtration of the process, and are intended to fill in the gap between the so-called structural and reduced form models for the default time [14]. In the credit risk context this particular model has the following interpretation [15]. $X$ represents the asset value of a company, and the default happens when $X$ reaches a default barrier, say, $x_1$. The market cannot observe $X$ perfectly, and is informed by the company management when there are significant changes in its economic standing, that is, when $X$ reaches some other thresholds, say, $x_2, \ldots, x_N$.

The purpose of this paper is to lay out the mathematics of this model. More precisely, we would like to describe the martingales and totally inaccessible stopping times of $\mathbb{F}^0$.

Our departure point is the Azéma's martingale. Azéma [1] established a formula for what is now known as the Azéma's martingale: let $B$ be a standard Brownian motion and define a filtration $\mathcal{A}_t^0 = \sigma\{\text{sign}(B_s); s \leq t\}$ and let $\mathbb{A}$ denote the completed filtration $(\mathcal{A}_t)_{t \geq 0}$. Let $A$ denote the martingale $A_t = E\{B_t | \mathcal{A}_t\}$ and letting $g_t = \sup\{s \leq t : B_s = 0\}$, Azéma's formula for this martingale is

$$A_t = \text{sign}(B_t) \sqrt{\frac{\pi}{2}} \sqrt{t - g_t}.$$

Emery [6] has named this Azéma's martingale and has proved many properties of it and related martingales, including homogeneous chaos representation and a fortiori martingale representation. (See [20] for an exposition and relevant background references.) There is, however, a treatment of a closely related theory established much earlier, related to regenerative sets in Markov process theory. Loosely speaking, a regenerative set is a homogeneous random set $M \subset \Omega \times [0, \infty)$ which has a renewal property at stopping times whose graphs are contained in $M$. For example, the zero set $\{(\omega, t) : B_t(\omega) = 0\}$ of Brownian motion is a regenerative set due to the strong Markov property. There is a classification of regenerative sets due to [12]. According to this classification, the zero set of Brownian motion is of a particular kind; for fixed $\omega$, $M(\omega)$ is a perfect and closed set with no interior. (We call these sets perfect regenerative sets.) We refer the reader to [12] for a precise definition, general classification of regenerative sets.

In [12], Jacod and Mémin focus on an unbounded perfect regenerative set $M \subset \Omega \times [0, \infty)$. They define the fundamental process $U_t(\omega) = t - \sup\{s \leq$ $t, (\omega, s) \in M\}$, and consider the random measure $\mu_0$ associated to the jumps of the process $U$. Taking $\mathbb{H} = (\mathcal{H}_t)_{t \geq 0}$ to be the smallest right continuous filtration which makes $U$ adapted, they provide the explicit form of the predictable compensator of $\mu_0$ with respect to $\mathbb{H}$ and prove that $\mu_0$ has the martingale representation property for $\mathbb{H}$; that is, any martingale $M$ adapted to $\mathbb{H}$ can be written as a stochastic integral with respect to $(\mu_0 - \nu_0)$. To use this to prove the martingale representation property of Azéma's martingale (by taking $M$ to be the zero set of the Brownian motion) requires minor modifications, since $\mathbb{H}$ does not give the same information as the previously defined filtration $\mathbb{A}$ (because one does not know when the paths are positive or negative, only when they are not zero). The key idea remains the same, however, that is, to make a time change of $\mathbb{A}$ using the inverse local time at zero, resulting in a filtration $\hat{\mathbb{A}}$ which can be shown to be generated by a Poisson random measure and, hence, to have martingale representation.

We build our results on the ideas of [12] because our setting can be thought of as a generalization of theirs [if we have only one level, say, $x_1$, then $\mathbb{F}^0$ is the same as the filtration generated by $(\text{sign}(X_s - x_1))_{s \geq 0}$]. Analogously, we prove a martingale representation theorem in terms of $N$ random measures, $\mu_1, \ldots, \mu_N$. These random measures are the key objects; we construct them and find their compensators using the technique of [12], [16] and [21], and the excursion theory of Markov processes. We also give formulae to explicitly compute the compensators in terms of the infinitesimal generator of the given diffusion. To do this, we use the results of [10] and [19].

Organization of this paper is as follows. In Section 2 we review the relevant results of excursion theory of Markov processes. In Section 3 we define the random measures $\mu_i, i = 1, \ldots, N$, and compute their compensators with respect to $\mathbb{F}$. In Section 4 we prove the representation theorem for $\mathbb{F}$, and in Section 5 we give a method to compute the compensators in terms of the infinitesimal generator of the underlying diffusion.

**2. Point process of excursions.** Here we give only a brief account of the excursion theory of Markov processes following [3]. We formulate the results in terms of the diffusion $X$ defined in the preceding section; however, we remark that this theory is developed for more general Markov processes (e.g., not necessarily with continuous paths). A detailed treatment of this subject can be found in [3]. Following the standard process notation, we consider the canonical six-tuple $(\Omega, \mathcal{G}, \mathbb{G}, X_t, \theta_t, P^x)$ describing the diffusion $X$.

Let $a \in \text{int}(I)$. We let $M_a(\omega) = \{t \geq 0 : X_t(\omega) = a\}$. $M_a(\omega)$ is a closed set since $X$ has continuous paths. As $(M_a(\omega))^c$ is open, there exist countably many disjoint open intervals, $(l_n, r_n)$, $n \geq 1$, such that

$$(M_a)^c = \bigcup_{n=1}^{\infty} (l_n, r_n).$$

(Here the dependence of all the terms on $\omega$ is implicit and dropped from the notation.) The open intervals $(l_n, r_n)$ are called the excursion intervals away from $a$. Note if $(l, r)$ is an excursion interval, the path $f(s) = X_{l+s}$, $0 \leq s < r - l$, is a continuous function whose graph lies entirely either in the upper positive half plane $\mathbb{R}_+ \times [a, \infty)$ or in the lower positive half plane $\mathbb{R}_+ \times (-\infty, a]$.

Let $L^a$ be the local time at $a$ and $\tau^a$ be its right continuous inverse, that is,
$$\tau_t^a = \inf\{s : L_s^a \geq t\}.$$
The key observation here is that $\Delta \tau_t^a > 0$ if and only if there exists an excursion interval $(l_n, r_n)$ such that $\tau_{t-}^a = l_n$. Using this, we can define a point process on $(\Omega, \mathcal{G}, P^a)$, called the excursion point process, denoted by $Y^a$. For fixed $\omega$, $Y^a(\omega)$ is the point function that maps the points in the set
$$D_{Y_\omega^a} = \{t > 0 | \tau_t^a(\omega) - \tau_{t-}^a(\omega) > 0\}$$
to $\Omega$ (i.e., the space of all continuous functions from $[0, \infty)$ to $I$). For each $t$ in $D_{Y_\omega^a}$, $Y_t^a(\omega)$ is defined as
$$Y_t^a(\omega)(s) = X^a(\theta_{\tau_{t-}^a}(\omega), s), \qquad s \in \mathbb{R}_+,$$
where $X^a$ denotes the stopped process $X$ at the hitting time, $H_a$, of $a$, that is,
$$X^a(\omega, t) = X(\omega, H_a(\omega) \wedge t),$$
where
$$H_a(\omega) = \inf\{t > 0, \omega(t) = a\}.$$
Note $Y_t^a$ is precisely the continuous path $f(s) = X_{l+s}$, $0 \leq s < r - l$, if $(l, r)$ is the excursion interval corresponding to $t$.

The following theorem, which is originally due to Itô [9], is the most important result of excursion theory. We give a generalized version of his result [17, 18] to include the case when $a$ is transient.

THEOREM 1 (Itô's theorem). *Under $P^a$, $Y^a$ is a Poisson point process absorbed at $\Omega_\infty = \{\omega \in \Omega | H_a(\omega) = \infty\}$ with respect to $\hat{\mathbb{G}}^a := (\mathcal{G}_{\tau_t^a})_{t \geq 0}$. Equivalently, there exists a measure $n^a$ on $(\Omega, \mathcal{G}^0)$ (called the characteristic measure of $Y^a$) s.t. for every nonnegative $\hat{\mathbb{G}}^a$-predictable process $Z$ and $(\Omega, \mathcal{G}^0)$-measurable nonnegative $g$*
$$E^a \sum_{\Delta \tau_s^a \neq 0} Z_s g(Y_s^a) = \int_\Omega g(\omega) n^a(d\omega) E^a \int_0^{L_\infty^a} Z_s \, ds.$$

We will refer to this formula as the *excursion formula*.

We go over two consequences of the excursion formula that are of particular importance to our analysis.

2.1. *Decomposition of the inverse local time and a representation theorem.*
Itô and McKean [10] gave the following decomposition for $\tau^a$ under $P^a$:

$$\tau_t^a = \tau_t^{a,-} + \tau_t^{a,+},$$

where

$$\tau_t^{a,+} = \int_0^{\tau_t^a} 1_{\{X_s > a\}} \, ds \quad \text{and} \quad \tau_t^{a,-} = \int_0^{\tau_t^a} 1_{\{X_s < a\}} \, ds.$$

When $a$ is recurrent, [10] showed that $\tau^{a,-}, \tau^{a,+}$ are two independent subordinators. A different formulation of this result can be given in terms of the point process of excursions which extends the result to the case when $a$ is transient. Note that the jump times of $\tau^a$ are exactly the jump times of $\tau^{a,-}$ and $\tau^{a,+}$; $\tau^{a,-}$ (resp. $\tau^{a,+}$) jumps at $t \in D_{Y_\omega^a}$ if and only if the graph of the excursion $Y_t^a$ lies entirely in $\mathbb{R}_+ \times (-\infty, a]$ (resp. $\mathbb{R}_+ \times [a, \infty)$). Let $U^{a,+} = \{e \in \Omega, e(t) \geq a, \forall t \in \mathbb{R}_+\}$ and $U^{a,-} = \{e \in \Omega, e(t) \leq a, \forall t \in \mathbb{R}_+\}$. Then $\Delta \tau_t^{a,\pm} = H_a(Y_t^a) 1_{\{Y_t^a \in U^{a,\pm}\}}$. So, we can consider a point process $\varphi^a$, defined for fixed $\omega$ as the point function from $D_{Y_\omega^a}$ to $[-\infty, 0) \cup (0, \infty]$, where

$$\varphi_t^a = H_a(Y_t^a) 1_{\{Y_t^a \in U^{a,+}\}} - H_a(Y_t^a) 1_{\{Y_t^a \in U^{a,-}\}} \qquad \text{for } t \in D_{Y_\omega^a}.$$

We let $\hat{\mathbb{H}}^a$ be the filtration generated by $\varphi^a$, made to be right continuous and complete with respect to $P^a$ (see [3] for a review of point processes).

Let $Z$ be $\hat{\mathbb{H}}^a$-predictable, and $g$ be Borel on $[-\infty, 0) \cup (0, \infty]$, both nonnegative. Since $Z$ is also $\hat{\mathbb{G}}^a$-predictable, the excursion formula gives

$$\begin{aligned}
(1) \quad & E^a \sum_{s \leq t} Z_s g(\varphi_s^a) 1_{\{\Delta \tau_s^a \neq 0\}} \\
& = E^a \int_0^{L_\infty^a \wedge t} Z_s \, ds \left( \int_{U^{a,+}} dn^a(e) g(H_a(e)) + \int_{U^{a,-}} dn^a(e) g(-H_a(e)) \right).
\end{aligned}$$

Let $F^a$ be the measure defined on $[-\infty, 0) \cup (0, \infty]$ by

$$F^a(A) = \int_{U^{a+}} dn^a(e) 1_{\{H_a(e) \in A \in (0, \infty]\}} + \int_{U^{a-}} dn^a(e) 1_{\{-H_a(e) \in A \in [-\infty, 0)\}}.$$

(1) implies that $\varphi^a$ is a Poisson point process absorbed at $\{-\infty\} \cup \{\infty\}$ with characteristic measure $F^a$.

Consider the random measure

$$\hat{\mu}^a(dt \times dx) = \sum_{\Delta \tau_s^a \neq 0} \varepsilon_{(s, \varphi_s^a)}(dt \times dx).$$

(1) also implies that the $\hat{\mathbb{H}}^a$-predictable compensator of $\hat{\mu}^a$ is the measure

$$\hat{\nu}^a(dt \times dx) = 1_{\{t \leq L_\infty^a\}} dt \times F^a(dx).$$

Since this uniquely characterizes the law of $\varphi^a$ (hence, $P^a$ on $\hat{\mathcal{H}}_\infty^a$) [18], $\hat{\mu}^a$ has the representation property (see, e.g., [13], page 174) for $\hat{\mathbb{H}}^a$, which we state in the following theorem.

Let $\mathcal{G}^1(\hat{\mu}^a)$ be the space of $\tilde{\mathcal{P}}_{\mathbb{H}}^a = \hat{\mathcal{P}}_{\mathbb{H}}^a \otimes \mathcal{B}([-\infty, 0) \cup (0, \infty])$-measurable functions defined on $\tilde{\Omega} = \Omega \times [0, \infty) \times [-\infty, 0) \cup (0, \infty]$ s.t. $E^a[(W^2 * \hat{\mu}_\infty^a)^{1/2}] < \infty$. (Here $\hat{\mathcal{P}}_{\mathbb{H}}^a$ is the predictable $\sigma$-algebra associated to $\hat{\mathbb{H}}^a$.)

THEOREM 2. *Any local martingale $\hat{M}$ adapted to $\hat{\mathbb{H}}^a$ is of the form*

$$\hat{M}_t = \hat{M}_0 + W * (\hat{\mu}^a - \hat{\nu}^a)_t,$$

*where $W \in \mathcal{G}_{\text{loc}}^1(\hat{\mu}^a)$.*

2.2. *Weil's formula for conditioning on $\hat{\mathcal{G}}_{T-}$.* Let $Z$ be a stochastic process and $\mathbb{A}$ be its natural filtration. Weil [21] gave a formula for the conditional distribution of $(Z_{T-}, Z_T)$ given the "strict past," $\mathcal{A}_{T-}$, for certain stopping times $T$ of $\mathbb{A}$, when $Z$ admits a Lévy system. [16] applied this formula to find the transition kernel of the process $U$ of a regenerative set $M$.

We reformulate Weil's original result for the point process of excursions of $X$. We omit the proof and refer the reader to [21] to see how easily his proof can be recast to prove our version using the excursion formula.

Let $\Phi$ be a nonnegative $\mathcal{B}[0, \infty) \otimes \mathcal{G}^0$-measurable function on $[0, \infty) \times \Omega$.

THEOREM 3 (Weil's formula). *Let $T = \inf\{s \in D_{Y_\omega^a} : (\tau_{s-}^a, Y_s^a) \in H\}$, where $H \in \mathcal{B}(0, \infty) \otimes \mathcal{G}^0$ is such that on $\{T < \infty\}$, $T(\omega) \in D_{Y_\omega^a}$ and $(\tau_{T-}^a, Y_T^a) \in H$ $P^a$-a.s. Then*

$$E^a[\Phi(\tau_{T-}^a, Y_T^a) 1_{\{0 < T < \infty\}} | \hat{\mathcal{G}}_{T-}^a] = Q_H(\tau_{T-}^a, \Phi) 1_{\{0 < T < \infty\}},$$

*where*

$$Q_H(t, \Phi) = \begin{cases} \dfrac{\int \Phi(t, e) 1_H(t, e) n^a(de)}{\int 1_H(t, e) dn^a(e)}, & \text{if } 0 < \int 1_H(t, e) dn^a(e) < \infty, \\ 0, & \text{if } \int 1_H(t, e) dn^a(e) = 0 \text{ or } \infty. \end{cases}$$

*Moreover, for any $H' \in \mathcal{B}(0, \infty) \otimes \mathcal{G}^0$, if $T_{H'}$ is defined as*

$$T_{H'}(\omega) = \begin{cases} T(\omega), & \text{if } T < \infty \text{ and } (\tau_{T-}^a, Y_T^a) \in H', \\ \infty, & \text{otherwise}, \end{cases}$$

*then*

$$E^a[\Phi(\tau_{T_{H'}-}^a, Y_{T_{H'}}^a) 1_{\{T_{H'} < \infty\}} | \hat{\mathcal{G}}_{T_{H'}-}^a] = Q_{H \cap H'}(\tau_{T_{H'}-}^a, \Phi) 1_{\{T_{H'} < \infty\}}.$$

**3. The filtration generated by $R(X)$ and the associated random measures.** Let us put $\mathcal{L} = \{x_1, \ldots, x_N\}$. Following [12], we can consider the process $U$ defined by

(2) $$U_t := t - \sup\{s \leq t : X_s \in \mathcal{L}\}.$$

If $\mathcal{L}$ consists of only one point, say, $a$, then $\mathbb{F}$ is the same as the filtration generated by $\text{sign}(X_s - a)$, and as discussed earlier in this case, the results of [12] can be modified to obtain a representation theorem for $\mathbb{F}$ with respect to the random measure associated to the jumps of $U$. Let us denote this random measure by $\mu$. When $\mathcal{L}$ consists of more than one point, $\mu$ is no longer enough to generate all martingales of $\mathbb{F}$. We can explain this using excursions. As before, $\{s \in [0, \infty) : X_s \in \mathcal{L}\}$ is a closed set, hence, its complement can be written as a countable union of disjoint open intervals $(a_n, b_n)$, we call these intervals excursion intervals away from $\mathcal{L}$. A representation theorem with respect $\mu$ implies the following: All martingales jump only at the ends of excursions, that is, at the points $b_n$ (since these are the only jump times of $U$), and the size of the jump (which could be zero) is determined by a predictable process. This is not true if $\mathcal{L}$ has more than one element because one can construct martingales that jump only if $X$ ends the excursion at a point $x_i$, the information of which is not "predictable." This suggests that we should decompose $\mu$ according to the values of $X$ at the $b_n$, which would give rise to random measures $\mu_i$ for each $x_i$. These random measures turn out to be the right choice in terms of which we can prove a martingale representation theorem.

We first establish our stochastic basis $(\Omega, \mathcal{F}, \mathbb{F}, P)$ and then construct the random measures $\mu_1, \ldots, \mu_N$ and find their compensators using the technique of [12] and [16].

We let $\mathcal{F}^i = (\bigvee_{t \geq 0} \mathcal{F}^0_t) \vee \mathcal{N}^i$, where $\mathcal{N}^i$ is the set of null sets of $P^{x_i}$. To have a complete stochastic basis for each $P^{x_i}$, we let $\mathcal{F}^i_t = \mathcal{F}^0_t \vee \mathcal{N}^i$.

PROPOSITION 4. $(\Omega, \mathcal{F}^i, \mathbb{F}^i, P^{x_i})$ *is a complete stochastic basis.*

PROOF. Clearly $\mathcal{F}^i$ is complete under $P^{x_i}$ and each $\mathcal{F}^i_t$ contains all $P^{x_i}$-null sets of $\mathcal{F}$. So we only need to check the right continuity of $\mathbb{F}^i$. First we observe that $\mathcal{F}^i_0 = \mathcal{F}^i_{0+}$, following from Blumenthal's 0–1 law, since $\mathcal{F}^0_{0+} \subset \mathcal{G}_{0+} \subset \mathcal{N}$. Next for $t > 0$, it is a standard result that $P^{x_i}(X_t \in \mathcal{L}) = 0$, therefore, $\{X_t \in \mathcal{L}\}$ is null, hence, in $\mathcal{F}^i_t$ and therefore, $(\bigcap_{s>t} \mathcal{F}^i_s) \cap \{X_t \in \mathcal{L}\} = \mathcal{F}^i_t \cap \{X_t \in \mathcal{L}\}$. For $t > 0$ s.t. $X_t \in \mathcal{L}^c$, $R(X)$ is constant in a neighborhood of $t$. In [4], page 304, this is referred to as strong continuity at $t$, by the argument of [4], page 304, we also have that

$$\left(\bigcap_{s>t} \mathcal{F}^i_s\right) \cap \{X_t \in \mathcal{L}^c\} = \mathcal{F}^i_t \cap \{X_t \in \mathcal{L}^c\}.$$

Now for any $A \in \bigcap_{s>t} \mathcal{F}_s^i$, both $A \cap \{X_t \in \mathcal{L}\} \in \mathcal{F}_t^i$ and $A \cap \{X_t \in \mathcal{L}^c\} \in \mathcal{F}_t^i$, hence, $A \in \mathcal{F}_t^i$. □

It will prove fundamental to have a collection of stopping times exhausting the jumps of $U$. This can be achieved as follows. For $x > 0$, let us define the following collection of stopping times:

$$T_1^x = \inf\{t > 0 : U_t > x\},$$
$$S_n^x = \inf\{t > T_n^x : U_t = 0\} \quad \text{for } n \geq 1,$$
$$T_{n+1}^x = \inf\{t > S_n^x : U_t > x\} \quad \text{for } n \geq 1.$$

Let $D = \{(\omega, t) : \Delta U_t(\omega) \neq 0\}$. Then $D = \bigcup_{x \in \mathbb{Q}} \bigcup_{n=1}^\infty [\![S_n^x]\!]$. We note that $U$ is adapted to $\mathbb{F}^i$ for each $i$, since $P^i$-a.s.

(3) $$U_t = t - \sup\{s \leq t, R(X_s) \neq R(X_t)\}.$$

We remark that (3) would not be true if the points $x_i$ were not regular, as this would not rule out the possibility that $R(X)$ is constant in a left neighborhood of $\sup\{s \leq t, X_s \in \mathcal{L}\}$.

As $U$ is adapted to $\mathbb{F}^i$, for each $n \geq 1$ and $x > 0$, $T_n^x$ and $S_n^x$ are stopping times of $\mathbb{F}^i$. We also observe that $X_{S_n^x}$ is $\mathcal{F}_{S_n^x}^i$-measurable, since

$$X_{S_n^x} = \sum_{i=1}^N x_{i+1} 1_{\{R(X_{S_n^x}) = R(X_{T_n^x}) = i\}} + x_i 1_{\{R(X_{S_n^x}) + 1 = R(X_{T_n^x}) = i\}}.$$

LEMMA 5.

(i) $\mathcal{F}_{T_n^x-}^i = \mathcal{F}_{T_n^x}^i$,
(ii) $\mathcal{F}_{S_n^x}^i = \mathcal{F}_{T_n^x}^i \vee \sigma(X_{S_n^x}, \Delta U_{S_n^x})$ (with the convention that $X_\infty = \infty$, and $\Delta U_\infty = \infty$),
(iii) $\mathcal{F}_{S_n^x-}^i = \mathcal{F}_{T_n^x}^i \vee \sigma(\Delta U_{S_n^x})$,
(iv) $\mathcal{F}_t^i \cap \{T_n^x \leq t < S_n^x\} = \mathcal{F}_{T_n^x}^i \cap \{T_n^x \leq t < S_n^x\}$.

PROOF. (i) This is an elementary consequence of the strong continuity of $R(X)$ at $T_n^x$. The details are omitted; see [4], T28, for a sample argument.

(ii) Since $X_{S_n^x}$ and $\Delta U_{S_n^x}$ are $\mathcal{F}_{S_n^x}^i$-measurable, it is enough to show $\mathcal{F}_{S_n^x}^i \subset \mathcal{F}_{T_n^x}^i \vee \sigma(X_{S_n^x}, \Delta U_{S_n^x})$. We observe $\mathcal{F}_{T_n^x}^i \vee \sigma(X_{S_n^x}, \Delta U_{S_n^x}) \supset \sigma(R(X_{s \wedge S_n^x}), s \geq 0) \vee \mathcal{N}^i$, since

$$R(X_{S_n^x \wedge t}) = R(X_{T_n^x \wedge t}) 1_{\{t < T_n^x - \Delta U_{S_n^x} - x\}}$$
$$+ 1_{\{t \geq T_n^x - \Delta U_{S_n^x} - x\}}$$
$$\times \left( \sum_i^N (i) 1_{\{R(X_{T_n^x}) = i, X_{S_n^x} = x_{i+1}\}} + (i-1) 1_{\{R(X_{T_n^x}) = i, X_{S_n^x} = x_i\}} \right).$$

So, we show
$$\sigma(R(X_{s \wedge S_n^x}), s \geq 0) \vee \mathcal{N}^i = \mathcal{F}_{S_n^x}^i.$$

In general, one has the following:
$$\mathcal{F}_{S_n^x}^i = \bigcap_{r>0} \sigma(R(X_{s \wedge (S_n^x+r)}), s \geq 0) \vee \mathcal{N}^i.$$

Note that
$$\bigcap_{r>0} \sigma(R(X_{s \wedge S_n^x+r}), s \geq 0) \subset \sigma(R(X_{s \wedge S_n^x}), s \geq 0) \vee \mathcal{G}_{0+}^{S_n^x},$$

where
$$\mathcal{G}_{0+}^{S_n^x} \doteq \sigma(\theta_{S_n^x}^{-1}(B) \cap \{S_n^x < \infty\} : B \in \mathcal{G}_{0+}).$$

We are done once we show
$$\mathcal{G}_{0+}^{S_n^x} \subset \sigma(X_{S_n^x}) \vee \mathcal{N}^i.$$

Let $A = \theta_{S_n^x}^{-1}(B) \cap \{S_n^x < \infty\}$. By the strong Markov property,
$$P^{x_i}(A \cap \{X_{S_n^x} = x_j\}) = P^{x_i}(X_{S_n^x} = x_j) P^{x_j}(B).$$

By Blumenthal's 0–1 law, $P^{x_j}(B) = 0$ or 1. If $P^{x_j}(B) = 0$, then $A \cap \{X_{S_n^x} = x_i\}$ is null, and if $P^{x_j}(B) = 1$, then $P^{x_i}(A \cap \{X_{S_n^x} = x_i\}) = P^{x_i}(X_{S_n^x} = x_i)$, implying that $\{X_{S_n^x} = x_i\} - (A \cap \{X_{S_n^x} = x_i\})$ is null. Therefore, in either case, $A \cap \{X_{S_n^x} = x_i\} \in \sigma(X_{S_n^x}) \vee \mathcal{N}^i$ and thus, $A \in \sigma(X_{S_n^x}) \vee \mathcal{N}^i$.

(iii) This is trivial because $\mathcal{F}_{S_n^x-}^i$ is generated by null sets and the sets of the form (for $s < t$)
$$\{R(X_s) = j\} \cap \{S_n^x > t\} = (\{R(X_s) = j\} \cap \{S_n^x > t\} \cap \{s \leq T_n^x\})$$
$$\cup (\{R(X_{T_n^x}) = j\} \cap \{S_n^x > t\} \cap \{T_n^x < s < S_n^x\}),$$

where the right-hand side is in $\mathcal{F}_{T_n^x} \vee \sigma(\Delta U_{S_n^x})$.

(iv) Trivial, because on $\{T_n^x \leq t < S_n^x\}$, $R(X_s) = R(X_{T_n^x \wedge s})$ for all $s \leq t$ □

Next, for each $i = 1, \ldots, N$, we define the stopping times
$$S_{n,i}^x = \begin{cases} S_n^x, & \text{if } X_{S_n^x} = x_i, \\ \infty, & \text{otherwise,} \end{cases}$$

and the corresponding random measures
$$\mu_i(\omega, ds, dx) = \sum_{s \geq 0} 1_{D_i}(\omega, s) \varepsilon_{(s, \Delta U_s)}(ds, dx)$$

on $\mathbb{R}_+ \times (0, \infty)$, where
$$D_i = \bigcup_{r \in \mathbb{Q}^+} \bigcup_{n=1}^{\infty} [|S_{n,i}^r|].$$

3.1. *Compensators of the random measures $\mu_i$.* For the rest of the discussion we fix $P = P^{x_{i^*}}$ and $\mathbb{F} = \mathbb{F}^{i^*}$. Let $\mathcal{P}$ be the predictable $\sigma$-algebra of $\mathbb{F}$. Our first task is to determine the predictable compensators of the $\mu_i$ under $P$, which we denote by $\nu_i$. Clearly, the process $\nu_i([0,t] \times (x,\infty))$ is the sum of its increments over the intervals $(T_n^x \wedge t, S_n^x \wedge t]$. Because of Lemma 5, these increments can be given in terms of the conditional distribution of $(\Delta U_{S_n^x}, X_{S_n^x})$ given $\mathcal{F}_{T_n^x}$ (see, e.g., [12, 13]). In particular,

$$\nu_i([0,t] \times (x,\infty))$$
$$= \begin{cases} 0, & \text{if } t \leq T_1^x, \\ \nu_i([0,T_n^x] \times (x,\infty)) + \int \frac{G_{n,x}^i(du)}{G_{n,x}((u,\infty])} 1_{\{x < u \leq x - T_n^x + t \wedge S_n^x\}}, \\ & \text{if } T_n^x \leq t \leq T_{n+1}^x, \end{cases}$$

where on $\{T_n^x < \infty\}$,

$$G_{n,x}^i(\omega, dy) = P(S_n^x - T_n^x + x \in dy, X_{S_n^x} = x_i | \mathcal{F}_{T_n^x})$$

and

$$G_{n,x}(\omega, dy) = \sum_{i=1}^N G_{n,x}^i(\omega, dy).$$

Note that, on $\{T_n^x < \infty\}$,

$$G_{n,x}(\omega, dy) = P(S_n^x - T_n^x + x \in dy | \mathcal{F}_{T_n^x}).$$

We use excursion theory to compute $G_{n,x}^i$, following the technique of [16]. For each $x_i$, let $U_i^+$ (resp. $U_i^-$) denote the set of all paths in $\Omega$ whose graphs lie in the upper (resp. lower) half plane $\mathbb{R}_+ \times [x_i, \infty)$ (resp. $\mathbb{R}_+ \times (-\infty, x_i]$). Let $n^i$ be the characteristic measure of the point process of excursions of $X$ away from $x_i$ and $H_i$ be the hitting time of $x_i$ [i.e., $H_i(\omega) = \inf\{t > 0 : \omega(t) = x_i\}$]. We define the following measures on $(0, \infty]$: For $A \in \mathcal{B}(0, \infty]$,

(4) $F_i^{1+}(A) = \int_{U_i^+} 1_{\{H_{i+1} \in A\}}(e) 1_{\{H_{i+1} < \infty\}}(e) \, dn^i(e) \quad$ for $i = 1, \ldots, N-1$,

(5) $F_i^{1-}(A) = \int_{U_i^-} 1_{\{H_{i-1} \in A\}}(e) 1_{\{H_{i-1} < \infty\}}(e) \, dn^i(e) \quad$ for $i = 2, \ldots, N$,

(6) $F_i^{0\pm}(A) = \int_{U_i^\pm} 1_{\{H_i \in A\}}(e) 1_{\{H_i \leq H_{i\pm 1}\}}(e) \, dn^i(e) \quad$ for $i = 1, \ldots, N$,

and

(7) $\qquad\qquad F_i^+(A) = F_i^{1+}(A) + F_i^{0+}(A),$

(8) $\qquad\qquad F_i^-(A) = F_i^{1-}(A) + F_i^{0-}(A),$

where we set $F_1^{1-} = F_N^{1+} = 0$.

Each of the measures $F_i^{j\pm}$ is a Lévy measure and the functions $F_i^{j\pm}[x, \infty]$ are continuous in $x$ (see Section 5).

THEOREM 6. $P$-a.s. on $\{T_n^x < \infty\}$,

$$G_{n,x}^i(\omega, dy) = \begin{cases} \dfrac{F_{i-1}^{1+}(dy)}{F_{i-1}^+(x, \infty]} 1_{\{y>x\}}, & \text{if } X_{T_n^x - x} = x_{i-1} \text{ and } R(X_{T_n^x}) = i - 1, \\ \dfrac{F_i^{0+}(dy)}{F_i^+(x, \infty]} 1_{\{y>x\}}, & \text{if } X_{T_n^x - x} = x_i \text{ and } R(X_{T_n^x}) = i, \\ \dfrac{F_i^{0-}(dy)}{F_i^-(x, \infty]} 1_{\{y>x\}}, & \text{if } X_{T_n^x - x} = x_i \text{ and } R(X_{T_n^x}) = i - 1, \\ \dfrac{F_{i+1}^{1-}(dy)}{F_{i+1}^-(x, \infty]} 1_{\{y>x\}}, & \text{if } X_{T_n^x - x} = x_{i+1} \text{ and } R(X_{T_n^x}) = i, \\ 0, & \text{otherwise.} \end{cases}$$

PROOF. We derive the first case only, other cases can be done similarly. We fix $n$ and $x$, and let $f$ be a nonnegative Borel function on $(0, \infty]$. If $T = \inf\{t > S_{n-1}^x, X_t = x_{i-1}\}$, by the strong Markov property,

$$
\begin{aligned}
&P(f(S_n^x - T_n^x + x) \cdot 1_{\{X_{S_n^x} = x_i\}} 1_{\{X_{T_n^x - x} = x_{i-1}, R(X_{T_n^x}) = i-1, T_n^x < \infty\}} | \mathcal{F}_{T_n^x}) \\
&(9) \quad = (P^{x_{i-1}}(f(S_1^x - T_1^x + x) 1_{\{X_{S_1^x} = x_i\}} 1_{\{X_{T_1^x - x} = x_{i-1}\}} 1_{\{R(X_{T_1^x}) = i-1\}} | \mathcal{F}_{T_1^x}^{i-1})) \\
&\quad \circ \theta_T 1_{\{T < T_n^x\}} 1_{\{T_n^x < \infty\}}.
\end{aligned}
$$

Let $\hat{T}_0 = \inf\{s \in D_{Y_\omega^{x_{i-1}}} : (\tau_{s-}^{x_{i-1}}, Y_s^{x_{i-1}}) \in [0, \infty) \times \Omega^x\}$, where $\Omega^x = \{\omega \in \Omega : H_i(\omega) \wedge H_{i-1}(\omega) \wedge H_{i-2}(\omega) > x\}$ and

$$\hat{T} = \begin{cases} \hat{T}_0, & \text{if } \hat{T}_0 < \infty, Y_{\hat{T}_0}^{x_{i-1}} \in U_{i-1}^+, \\ \infty, & \text{otherwise.} \end{cases}$$

Let $A = \{X_{T_1^x - x} = x_{i-1}, R(X_{T_1^x}) = i - 1\}$. We note that $1_{\{\hat{T} < \infty\}} = 1_A$, and on $A$, $T_1^x = \tau_{\hat{T}-}^{x_{i-1}} + x$, $P^{x_{i-1}}$ a.s. Moreover, on $A$, $R(X_u) = i$ if $T_1^x - x \leq u \leq T_1^x$, and $R(X_u) 1_{\{u < T_1^x - x\}}$ is $\hat{\mathcal{G}}_{\hat{T}-}^{x_{i-1}}$-measurable (since $u < T_1^x - x$ if and only if $L_u^{x_{i-1}} < \hat{T}$). Thus,

$$\mathcal{F}_{T_1^x}^{i-1} \cap A \subset \hat{\mathcal{G}}_{\hat{T}-}^{x_{i-1}} \vee \mathcal{N}^{i-1}.$$

Therefore, $P^{x_{i-1}}$-a.s. on $A$,

$$
\begin{aligned}
&P^{x_{i-1}}(f(S_1^x - T_1^x + x) 1_{\{X_{S_1^x} = x_i\}} 1_{\{X_{T_1^x - x} = x_{i-1}\}} | \mathcal{F}_{T_1^x}^{i-1}) \\
&(10) \quad = P^{x_{i-1}}(P^{x_{i-1}}(f(H_i(Y_{\hat{T}}^{x_{i-1}})) 1_{\{H_i(Y_{\hat{T}}^{x_{i-1}}) < \infty\}} 1_{\{\hat{T} < \infty\}} | \hat{\mathcal{G}}_{\hat{T}-}^{x_{i-1}}) | \mathcal{F}_{T_1^x}^{i-1}).
\end{aligned}
$$

By Weil's formula, the inner conditional expectation is given by

$$\tag{11} \frac{\int_{U_{i-1}^+} f \circ H_i 1_{\{x < H_i < \infty\}} \, dn^{i-1}}{\int_{U_{i-1}^+} 1_{\{H_i \wedge H_{i-1} > x\}} \, dn^{i-1}} 1_{\{\hat{T} < \infty\}} = \frac{\int_{(x,\infty]} f \, dF_{i-1}^{1+}}{\int_{(x,\infty]} dF_{i-1}^+} 1_{\{\hat{T} < \infty\}}.$$

The desired expression follows from (9), (10) and (11), and observing $1_{\{\hat{T} < \infty\}} \circ \theta_T 1_{\{T < T_n^x\}} 1_{\{T_n^x < \infty\}} = 1_{\{X_{T_n^x -} = x_{i-1}\}} 1_{\{X_{T_n^x} > x_{i-1}\}} 1_{\{T_n^x < \infty\}}.$ □

COROLLARY 7.

  (i) $P$-a.s., $\nu_i(\omega, \{t\} \times (0,\infty)) = 0$.
  (ii) The stopping times $S_n^x$ are totally inaccessible with respect to $\mathbb{F}$.

PROOF. (i) is true simply because the functions $F_i^{j\pm}[x,\infty]$ are continuous. For (ii), observe that the compensator of the one point process $1_{\{S_n^x \leq t\}}$ is

$$A_t = 1_{\{t \geq T_n^x\}} \sum_{i=1}^N \int_x^{(t \wedge S_n^x) - T_n^x + x} \frac{G_{n,x}^i(du)}{G_n^x[u,\infty]},$$

which is again continuous because the functions $F_i^{j\pm}[x,\infty]$ are continuous. □

**4. Representation theorem.** Let $\mathcal{M}$ (resp. $\mathcal{M}^2$) be the space of uniformly integrable (resp. square integrable) martingales with càdlàg paths and adapted to $\mathbb{F}$. Let $\mathcal{V}$ be the space of processes with finite variation, again with càdlàg paths and adapted to $\mathbb{F}$. Let $\mathcal{G}^{1,i}$ (resp. $\mathcal{S}_i^2$, $\mathcal{S}_i$) be the space of $\tilde{\mathcal{P}} := \mathcal{P} \otimes \mathcal{B}(0,\infty)$- measurable functions $W$ on $\tilde{\Omega} = \Omega \times [0,\infty) \times (0,\infty)$ s.t. $E[(W^2 * \mu_i)_\infty)^{1/2}] < \infty$ (resp. $E[(W^2 * \mu_i)_\infty] < \infty$, $E[(|W| * \mu_i)_\infty] < \infty$). (Note that, due to Corollary 7, these are the right choices of integrand spaces in order for the stochastic integral to be in $\mathcal{M}$, $\mathcal{M}^2$ and $\mathcal{M} \cap \mathcal{V}$, resp.; see, e.g., [13] for a definition of the stochastic integral with respect to a random measure.)

THEOREM 8. For any $M \in \mathcal{M}_{\mathrm{loc}}$, there exists $W_i$ in $\mathcal{G}_{\mathrm{loc}}^{1,i}$ such that

$$M_t = E[M_0] + \sum_{i=1}^N (W_i * (\mu_i - \nu_i))_t.$$

Moreover, if $M \in \mathcal{M}^2$, then $W_i \in \mathcal{S}_i^2$; and if $M \in \mathcal{M}$ and of finite variation, then $W_i \in \mathcal{S}_i$.

PROOF. Let us consider the random measure $\bar{\mu}$ on $\overline{\Omega} \doteq \Omega \times [0, \infty) \times E$, where $E = (0, \infty) \times \{x_1, \ldots, x_n\}$ defined by
$$\bar{\mu}(\omega, dt \times (dx \times \{x_i\})) = \mu_i(\omega, dt \times dx).$$
Let $\overline{\mathcal{P}} \doteq \mathcal{P} \times \mathcal{E}$, where $\mathcal{E}$ is the $\sigma$-algebra generated by the sets $B \times \{x_i\}$, and $B \in \mathcal{B}(0, \infty)$; and $\mathcal{P}$ is the predictable $\sigma$-algebra of $\mathbb{F}$. We define the spaces $\mathcal{G}^1$, $\mathcal{S}^1$ and $\mathcal{S}^2$ for $\bar{\mu}$ analogously. Clearly, $\overline{\mathcal{P}}$-compensator of $\bar{\mu}$ can be given in terms of $\nu_i$, in particular,
$$\nu(\omega, dt \times (dx \times \{x_i\})) = \nu_i(\omega, dt \times dx).$$
Let $D = \{(\omega, t) : \mu(\omega, \{t\} \times E) = 1\}$. Then we have that
$$D = \bigcup_{r \in \mathbb{Q}^+} \bigcup_{n=1}^{\infty} [|S_n^r|].$$
We recall that
$$(12) \qquad \mathcal{F}_{S_n^x} = \mathcal{F}_{T_n^x} \vee \sigma(X_{S_n^x}, \Delta U_{S_n^x}).$$
Since $\{X_{S_n^x} = x_i, -\Delta U_{S_n^x} \in A\} = \{\bar{\mu}([S_n^r] \times (A \times \{x_i\})) = 1\}$ and
$$D = \bigcup_{n,m} [|S_n^{1/m}|],$$
condition (12) allows us to use a special case of Jacod's decomposition theorem of martingales with respect to random measures (see [11], e.g., Theorem 4.1 and remark 3 following Theorem 4.1, and Proposition 2.4). That is, for any local martingale $M$, there exists $W \in \mathcal{G}_{\text{loc}}^1$ such that
$$M = M_0 + W * (\bar{\mu} - \bar{\nu}) + M',$$
where $M_0$ is $\mathcal{F}_0$-measurable and integrable, and $M'$ is a local martingale with $M_0' = 0$, adapted to $\mathbb{F}$, and that does not jump on the support of $\bar{\mu}$. Since $\mathcal{F}_0$ consists only of null sets, $M_0 = E[M_0]$ a.s. Also, if we let $W_i(\omega, t, y) = W(\omega, t, (y, x_i))$, then each $W_i \in \mathcal{G}_{\text{loc}}^{1,i}$ and
$$W * (\bar{\mu} - \bar{\nu}) = \sum_{i=1}^{N} W_i * (\mu_i - \nu_i).$$
According to [11], if $M \in \mathcal{M}$ and of finite variation (resp. $M \in \mathcal{M}^2$), then $W \in \mathcal{S}^1$ (resp. $\mathcal{S}^2$), which implies that each $W_i \in \mathcal{S}_i$ (resp. $\mathcal{S}_i^2$). So we only need to show that $M' = 0$ and by localization, we may assume that $M \in \mathcal{M}$.

We first argue that $M'$ has a constant stretch during any excursion interval. Equivalently, we can show $M'_{S_n^x} - M'_{T_n^x} = 0$ for any $x > 0$ and $n \geq 1$. We consider the filtration $\mathbb{F}^{n,x} = (\mathcal{F}_t^{n,x})_{t \geq 0}$, where $\mathcal{F}_t^{n,x} = \mathcal{F}_{T_n^x} \vee \mathcal{F}_{S_n^x \wedge t}$, and the multivariate "one" point process $\mu^{n,x}$ defined by
$$\mu^{n,x}(\omega, [0, t] \times \{x_i\}) = 1_{\{S_n^x - T_n^x + x \leq t, X_{S_n^x} = x_i\}}.$$

Due to Lemma 5, $\mathbb{F}^{n,x}$ is the smallest right continuous filtration $\mathbb{A}$ for which $\mu^{n,x}$ is optional and that has $\mathcal{A}_0 = \mathcal{F}_0^{n,x}$. Therefore, $\mathbb{F}^{n,x}$ has martingale representation with respect to $\mu^{n,x}$ (see [13], page 148). The martingale $M_t^{n,x} = M'_{(T_n^x+t)\wedge S_n^x}$ is adapted to $\mathbb{F}^{n,x}$ and does not jump at $S_n^x - T_n^x + x$, therefore has to be constant.

The second and last step is proving the following claim: Any uniformly integrable martingale that stays constant during the excursion intervals is constant on $[0,\infty)$. To prove this, we use a time change argument. To motivate, suppose $\mathcal{L}$ consists of only one point, $x_{i^*}$. The uniformly integrable martingale $\hat{M}'$ obtained by applying the time change $\hat{M}'_t = M'_{\tau_t}$, where $\tau$ is the inverse local time of $X$ at $x_{i*}$, does not jump on the support of the random measure $\hat{\mu}$ corresponding to the jumps of $\tau^+ - \tau^-$. Since $\hat{M}'$ is adapted to the filtration generated by $\hat{\mu}$, and $\hat{\mu}$ has the martingale representation property by Theorem 2, it follows that $\hat{M}'_t = \hat{M}'_0 = 0$. Thus, $M' = 0$, since $M'_t = E[\hat{M}'_{L_t}|\mathcal{F}_t]$. Now we improve this argument to prove the claim for the case when $\mathcal{L}$ consists of more than one point. Let $g_t = \sup\{s < t, X_s \in \mathcal{L}\}$ and consider the process $\hat{X}$ defined by $\hat{X}_t = X_{g_t}$. The paths of $\hat{X}$ are step functions taking values in $\mathcal{L}$. If $\hat{X}_t = x_i$, the next jump time of $\hat{X}$ is the first time $s > t$, such that $X_s \in [x_{i-1}, x_{i+1}]^c$. Therefore, if we consider the stopping times

$$T_0 = 0,$$
$$T_n = \inf\{t > T_{n-1} : \Delta \hat{X}_t \neq 0\}, \qquad n \geq 1,$$

we have $X_{T_n} \in \mathcal{L}$, and $\bigcup_{n=0}^{\infty}[T_n, T_{n+1}) = [0,\infty)$ due to continuity of $X$. We let

$$X_t^n = (X_{T_{n-1}+t\wedge T_n})1_{\{T_{n-1}<\infty\}}.$$

By the strong Markov property, $X^n$ is a copy of the original diffusion with the starting point $X_{T_{n-1}}$ stopped at the exit of the interval $[x_{i-1}, x_{i+1}]$ if $X_{T_{n-1}} = x_i$. Let $L^{n,i}$ be the local time of $X^n$ at $x_i$, $\tau^{n,i}$ be its right continuous inverse, and $\varphi^{n,i}$ be the point process of the jumps of $\tau^{n,i,+} - \tau^{n,i,-}$, where

$$\tau_t^{n,i,+} = \int_0^{\tau_t^{n,i}} 1_{\{X_s^n > x_i\}} ds \quad \text{and} \quad \tau_t^{n,i,-} = \int_0^{\tau_t^{n,i}} 1_{\{X_s^n < x_i\}} ds.$$

Let

$$L^n = \sum_{i=1}^{N} L^{n,i} 1_{\{X_{T_{n-1}}=x_i\}},$$

$$\tau^n = \sum_{i=1}^{N} \tau^{n,i} 1_{\{X_{T_{n-1}}=x_i\}},$$

$$\varphi^n = \sum_{i=1}^{N} 1_{\{X_{T_{n-1}}=x_i\}} \varphi^{n,i}$$

and we let $\mu^n$ be the random measure corresponding to the jumps of $\varphi^n$. We consider the filtration $\mathbb{H}^n = (\mathcal{H}^n_t)_{t \geq 0}$, where

$$\mathcal{H}^n_t = \sigma(\varphi^n_s, s \leq t) \vee \mathcal{F}_{T_{n-1}}.$$

Let $\mathcal{P}^n$ be the predictable $\sigma$-algebra associated to $\mathbb{H}^n$. The predictable compensator $\nu^n$ (with respect to $\mathcal{P}^n$) of $\mu^n$ is given by

$$\nu^n(dt \times dx) = \sum_{i=1}^N 1_{\{X_{T_{n-1}} = x_i\}} 1_{\{t \leq L^{n,i}_\infty\}} dt \times n^i_{X^n}(dx).$$

By the discussion preceding Theorem 2, it follows that $\nu^n$ completely describes the restriction of $P(\cdot | X_{T_{n-1}} = x_i)$ to $\sigma(\varphi^n_s, s \geq 0)$ for each $x_i$. Therefore, $P$ on $\mathcal{H}^n_\infty$ is determined by $\nu^n$ and the restriction of $P$ to $\mathcal{F}_{T_{n-1}}$, that is, following the notation of [13], the martingale problem $s(\mathcal{F}_{T_{n-1}}, \mu^N | P_{|\mathcal{F}_{T_{n-1}}}, \nu^n)$ has a unique solution. By the fundamental representation theorem (see, e.g., [13], page 212), this implies martingale representation with respect to $\mu^n$ for $\mathbb{H}^n$.

We let

$$N_t = (M'_{(T_{n-1}+t) \wedge T_n} - M'_{T_{n-1}}) 1_{\{T_{n-1} < \infty\}}$$

and

$$\hat{N}_t = \sum_{i=1}^N N_{\tau^{n,i}_t} 1_{\{X_{T_{n-1}} = x_i\}}.$$

Since $M'$ stays constant during the excursion intervals, we have that $\hat{N}_t - \hat{N}_{t-} = 0$ for all $t$ such that $\Delta \tau^n_t \neq 0$. Next we argue that $\hat{N}$ is adapted to $\mathbb{H}^n$. Let $\mathcal{F}^n_t = \mathcal{F}_{T_{n-1}+t}$. Then $\hat{N}_t$ is $\mathcal{F}^n_{\tau^n_t}$-measurable and $\mathcal{F}^n_{\tau^n_t} \cap \{t < L^n_\infty\} \subset \mathcal{H}^n_t \cap \{t < L^n_\infty\}$. The latter follows from the fact that $\mathcal{F}^n_t \cap \{L^n_t < L^n_\infty\} \subset \mathcal{H}_{L^n_t} \cap \{L^n_t < L^n_\infty\}$ since on $\{L^n_t < L^n_\infty\}$, one has $R(X_t) = R(X_{T_{n-1}})$ or $R(X_{T_{n-1}}) - 1$, which can be determined by the sign of $\varphi^n_{L^n_t}$. This implies that $\hat{N}_t 1_{\{t < L^n_\infty\}}$ is $\mathcal{H}^n_t$-measurable and $\hat{N}_{L^n_\infty -}$ is $\mathcal{H}^n_\infty$-measurable. Since $\hat{N}_{L_\infty -} = \hat{N}_{L_\infty}$, $\hat{N}$ is adapted to $\mathbb{H}^n$.

Adapted to $\mathbb{H}^n$ and having no jumps on the support of $\mu^n$, $\hat{N}$ must be equal to $\hat{N}_0(=0)$. Since $N_t = E[\hat{N}_{L^n_t} | \mathcal{F}^n_t]$, $N = 0$ and therefore, $M'_{T_n} - M'_{T_{n-1}} = 0$. Observing $\bigcup_{n=1}^\infty [T_{n-1}, T_n] = [0, \infty)$, it follows that $M'$ is constant (hence, is equal to $M'_0 = 0$). $\square$

COROLLARY 9.

(i) *All local martingales adapted to $\mathbb{F}$ are purely discontinuous with jumps only at the ends of excursions.*

(ii) *The graph $[|T|]$ of any totally inaccessible stopping time $T$ is contained in $D$.*

**5. Computation of the measures $F_i^{j\pm}$.** As we have seen, the compensators $\nu_i$ of the random measures $\mu_i$ are given in terms of the measures $F_i^{j\pm}$, which are, therefore, the key quantities that we want to compute.

The measures $F_i^{j\pm}$ are related to the Lévy measure of the inverse local time of $X$ at $x_i$. That is the measure $F_i$ defined on $(0,\infty]$ as $F_i(A) = \int_\Omega 1_{\{H_i \in A\}} dn^i$ (see, e.g., [19], Sections 2.1 and 2.2). Clearly, $F_i^{j\pm}[x,\infty] \leq F_i[x,\infty]$ $\forall x > 0$, hence, each $F_i^{j\pm}$ defines a Lévy measure on $(0,\infty]$ as well.

We let
$$\psi_i^{j\pm}(\lambda) = \int_{(0,\infty]} (1 - e^{-\lambda x}) F_i^{j\pm}(dx) = \lambda \int_0^\infty e^{-\lambda x} F_i^{j\pm}[x,\infty]\, dx.$$

In this section we prove a theorem (which we obtain as a corollary to the results of [19] on the decomposition of the Laplace exponent of the inverse local time and last exit times) characterizing $\psi_i^{j\pm}$ in terms of the solutions of the equation $\mathcal{A}\Phi = \lambda\Phi$.

We let $\Phi_{i,\lambda}^+$ for $i = 1, \ldots, N-1$ (resp. $\Phi_{i,\lambda}^-$ for $i = 2, \ldots, N$) denote the decreasing (resp. increasing) solution of $\mathcal{A}\Phi = \lambda\Phi$ with boundary condition $\Phi(x_{i+1}) = 0$ [resp. $\Phi(x_{i-1}) = 0$] normalized to have $\Phi_{i,\lambda}^\pm(x_i) = 1$. We also let $\tilde{\Phi}_{i,\lambda}^+$ (resp. $\tilde{\Phi}_{i,\lambda}^-$) denote the decreasing (resp. increasing) solution of $\mathcal{A}\Phi = \lambda\Phi$ with boundary condition $\Phi(\sup(I)) = 0$ [resp. $\Phi(\inf(I)) = 0$]. (Here we only deal with absorbing boundary conditions at the boundary of $I$.) (For the construction of such solutions, see [10], page 128.) For $\lambda > 0$, let
$$\psi_i^\pm(\lambda) = -\pm \frac{1}{2} \frac{(\Phi_{i,\lambda}^\pm)'(x_i)}{\Phi_{i,\lambda}^\pm(x_i)},$$

and similarly,
$$\tilde{\psi}_i^\pm(\lambda) = -\pm \frac{1}{2} \frac{(\tilde{\Phi}_{i,\lambda}^\pm)'(x_i)}{\tilde{\Phi}_{i,\lambda}^\pm(x_i)}.$$

We put $\psi_i^\pm(0) = \lim_{\lambda \to 0} \psi_i^\pm(\lambda)$ and $\tilde{\psi}_i^\pm(0) = \lim_{\lambda \to 0} \tilde{\psi}_i^\pm(\lambda)$. Finally we let
$$\psi_{i,i+1}(\lambda) = \psi_i^+(\lambda) + \tilde{\psi}_i^-(\lambda)$$

and
$$\psi_{i,i-1}(\lambda) = \psi_i^-(\lambda) + \tilde{\psi}_i^+(\lambda).$$

THEOREM 10. (i) *For $i = 1, \ldots, N-1$,*
$$\psi_i^{1+}(\lambda) = P_{x_i}(H_{i+1} < \infty)\psi_{i,i+1}(0) - P_{x_i}(e^{-\lambda H_{i+1}})\psi_{i,i+1}(\lambda),$$
*and for $i = 2, \ldots, N$,*
$$\psi_i^{1-}(\lambda) = P_{x_i}(H_{i-1} < \infty)\psi_{i,i-1}(0) - P_{x_i}(e^{-\lambda H_{i-1}})\psi_{i,i-1}(\lambda).$$

(ii) *For $i = 1, \ldots, N - 1$,*
$$\psi_i^{0+}(\lambda) = -\psi_i^+(0) + \psi_i^+(\lambda),$$
*and for $i = 2, \ldots, N$,*
$$\psi_i^{0-}(\lambda) = -\psi_i^-(0) + \psi_i^-(\lambda).$$
(iii) $\psi_N^{0+}(\lambda) = \tilde{\psi}_N^+(\lambda)$ *and* $\psi_1^{0-}(\lambda) = \tilde{\psi}_1^-(\lambda)$.

PROOF. (iii) directly follows from [19] (see, in particular, Theorem 1 and Section 2.2).

(ii) For each $i = 1, \ldots, N-1$, consider the diffusion $X^{i+1}$ with state space $I \cap (-\infty, x_{i+1}]$ obtained by stopping $X$ at $H_{i+1}$. The infinitesimal generator of $X^{i+1}$ coincides with $\mathcal{A}$ on the interior of $I \cap (-\infty, x_{i+1}]$ and $x_{i+1}$ becomes an absorbing boundary point. Let $n_{X^{i+1}}^i$ be the characteristic measure of the point process of excursions of $X^{i+1}$ around $x_i$. Then for $A \in \mathcal{B}(0, \infty]$,
$$F_i^{0+}(A) = \int_{U_i^+} 1_{\{H_i \in A\}}(e) 1_{\{H_i < \infty\}}(e) \, dn_{X^{i+1}}^i(e).$$
Let $F_i^+$ be the measure on $(0, \infty]$ defined by $F_i^+(A) = \int_{U_i^+} 1_{\{H_i \in A\}}(e) \, dn_i^+(e)$. According to Pitman and Yor [19] (see also [10], page 128 with regards to the boundary conditions), we have that
$$\psi_i^+(\lambda) = \int_{(0, \infty]} (1 - e^{-\lambda x}) F_i^+(dx).$$
Therefore, $\psi_i^{0+}(\lambda) = \psi_i^+(0) - \psi_i^+(\lambda)$, for $i = 1, \ldots, N-1$. A similar argument (replacing $+$ with $-$) yields the result for $\psi_i^{0-}$, as well, for $i = 2, \ldots, N$.

(i) For $\psi_i^{1+}$ for $i = 1, \ldots, N-1$ (resp. $\psi_i^{1-}$ for $i = 2, \ldots, N$), we are again going to make use of the results of [19]. Following their notation, let
$$\Lambda_{i,i+1} = \sup\{t < H_{i+1} : X_t = x_i\}$$
be the last exit time from $x_i$ before the first hit of $x_{i+1}$. Pitman and Yor [19] finds the Laplace transform of $H_{i+1} - \Lambda_{i,i+1}$ as
$$E^{x_i}(e^{-\lambda(H_{i+1} - \Lambda_{i,i+1})}) = E^{x_i}(e^{-\lambda H_{i+1}}) \frac{\psi_{i,i+1}(\lambda)}{\psi_{i,i+1}(0)}.$$
We observe that
$$E^{x_i}(1_{\{H_{i+1} < \infty\}} - e^{-\lambda(H_{i+1} - \Lambda_{i,i+1})})$$
$$= E^{x_i} \sum_{s \leq L_{H_{i+1}}^{x_i}, \tau_{s-}^{x_i} \neq \tau_s^{x_i}} (1 - e^{-\lambda H_{i+1}(Y^{x_i}(s))}) 1_{\{H_{i+1}(Y^{x_i}(s)) < \infty\}}$$
$$= E^{x_i}(L_{H_{i+1}}) \int (1 - e^{-\lambda H_{i+1}(e)}) 1_{\{H_{i+1} < \infty\}} n^i(de),$$

where the last equality follows from the excursion formula. We note that the second term in the right-hand side is $\psi_i^{1+}(\lambda)$. Since $L_{H_{i+1}}$ is equal to $L_\infty^{i+1}$ where $L^{i+1}$ is the local time at $x_i$ of the stopped diffusion $X^{i+1}$, from Theorem 2 of [19], we have that

$$E^{x_i}(L_{H_{i+1}}) = \frac{1}{\psi_{i,i+1}(0)}.$$

Now the result follows by simple algebra for $\psi_i^{1+}$ for $i = 1, \ldots, N-1$. By the same argument replacing $+$ with $-$ we also get the expression for $\psi_i^{1-}$ for $i = 2, \ldots, N$. $\square$

This theorem can be used to compute the measures $F_i^{j\pm}$, if the coefficients of $\mathcal{A}$ are parameterized, and if one can identify the solutions of $\mathcal{A}\Phi = \lambda\Phi$ for the given parametrization. However, what remains open is the qualitative understanding of how these measures react to the changes in the coefficients of the diffusion in a general setting where the coefficients are nonparameterized.

5.1. *Remarks.* In Section 3 we stated that the functions $F_i^{0,\pm}[x, \infty]$ and $F_i^{1\pm}[x, \infty)$ are continuous. Here we justify why this is correct.

Pitman and Yor [19] point out that (with reference to [2]) the measures $F_i^{0\pm}$ have spectral representation, that is, there exist measures $\mu_i^{0\pm}$ on $[0, \infty)$ with $\int_{[0,\infty)} (1+\xi)^{-1} \mu_i^\pm(d\xi) < \infty$ such that

$$F_i^{0\pm}[x, \infty] = \int_0^\infty e^{-\xi x} \mu_i^{0\pm}(d\xi).$$

It follows that $F_i^{0\pm}[x, \infty]$ is continuous, in fact, $F_i^{0\pm}$ has a smooth density with respect to Lebesgue the measure and the density is given by the function

$$f_i^{0\pm}(x) = \int_0^\infty \xi e^{-\xi x} \mu_i^{0\pm}(d\xi).$$

From the Laplace transform of $F_i^{1\pm}[x, \infty)$, we see that

$$F_i^{1\pm}[x, \infty) = P^{x_i}(H_{i\pm 1} < \infty) F_i^{i\pm 1}(\infty) - \int_0^x g_i^{i+1}(x-u)(F_i^{i+1}[u, \infty])\, du,$$

where $F_i^{i+1}$ is the Lévy measure of the inverse local time at $x_i$ of the stopped diffusion $X^{i+1}$ and $g_i^{i+1}$ is the probability density function of the hitting time of $x_{i+1}$ with respect to $P^{x_i}$, which exists and continuous is (see, e.g., [10]). Since $F_i^{i+1}[x, \infty]$ is also continuous by the same comment of [19] as above, we have the continuity of $F_i^{1\pm}[x, \infty)$.

**Acknowledgment.** I would like to thank my Ph.D. advisor Professor Philip Protter. This work has been completed under his thesis direction during my Ph.D. studies at the School of ORIE, Cornell University, Ithaca, NY.

DEPARTMENT OF MATHEMATICS AND STATISTICS
YORK UNIVERSITY
4700 KEELE ST.
TORONTO, ONTARIO
CANADA M3J 1P3
E-MAIL: dsezer@mathstat.yorku.ca